\documentclass[11pt]{amsart}
\usepackage{amsmath,amsfonts,amsthm,amsopn,amssymb}
\usepackage{cite,marginnote}
\usepackage{color,enumitem,graphicx}
\usepackage[colorlinks=true,urlcolor=blue,
citecolor=red,linkcolor=blue,linktocpage,pdfpagelabels,
bookmarksnumbered,bookmarksopen]{hyperref}
\usepackage[english]{babel}

\usepackage[left=2.7cm,right=2.7cm,top=2.9cm,bottom=2.9cm]{geometry}
\usepackage[hyperpageref]{backref}

\numberwithin{equation}{section}

\makeindex
\newtheorem{theorem}{Theorem}[section]
\newtheorem{proposition}[theorem]{Proposition}
\newtheorem{lemma}[theorem]{Lemma}
\newtheorem{remark}[theorem]{Remark}
\newtheorem{example}[theorem]{Example}
\newtheorem{corollary}[theorem]{Corollary}
\newtheorem{definition}[theorem]{Definition}

\newcommand{\ud}{\mathrm{d}}
\newcommand{\RN}{\mathbb R^N}

\newcommand{\iy}{\infty}

\newcommand{\s}{\section}

\newcommand{\DD}{\Delta}
\newcommand{\g}{\gamma}
\newcommand{\G}{\Gamma}

\newcommand{\la}{\lambda}

\newcommand{\si}{\sigma}

\newcommand{\R}{\mathbb R}
\newcommand{\al}{\alpha}

\newcommand{\ti}{\tilde}

\newcommand{\re}[1]{(\ref{#1})}

\newcommand{\rg}{\rightarrow}

\newcommand{\e}{\varepsilon}

\newcommand{\lab}{\label}
\newcommand{\bt}{\begin{theorem}}
\newcommand{\et}{\end{theorem}}
\newcommand{\bl}{\begin{lemma}}
\newcommand{\el}{\end{lemma}}
\newcommand{\bd}{\begin{definition}}
\newcommand{\ed}{\end{definition}}
\newcommand{\bc}{\begin{corollary}}
\newcommand{\ec}{\end{corollary}}
\newcommand{\bp}{\begin{proof}}
\newcommand{\ep}{\end{proof}}
\newcommand{\bx}{\begin{example}}
\newcommand{\ex}{\end{example}}
\newcommand{\bi}{\begin{exercise}}
\newcommand{\ei}{\end{exercise}}
\newcommand{\bo}{\begin{proposition}}
\newcommand{\eo}{\end{proposition}}
\newcommand{\br}{\begin{remark}}
\newcommand{\er}{\end{remark}}
\newcommand{\be}{\begin{equation}}
\newcommand{\ee}{\end{equation}}
\newcommand{\ba}{\begin{align}}
\newcommand{\ea}{\end{align}}
\newcommand{\bn}{\begin{enumerate}}
\newcommand{\en}{\end{enumerate}}
\newcommand{\bg}{\begin{align*}}
\newcommand{\bcs}{\begin{cases}}
\newcommand{\ecs}{\end{cases}}

\newcommand{\bean}{\begin{eqnarray*}}
\newcommand{\eean}{\end{eqnarray*}}

\title[Fractional Schr\"odinger-Poisson systems]{Fractional Schr\"odinger-Poisson systems with a \\ general
	subcritical or critical nonlinearity}

\author[J. Zhang]{Jianjun Zhang}
\author[J.M.\ do \'O]{Jo\~ao Marcos do \'O}
\author[M. Squassina]{Marco Squassina}

\address[J. J.\ Zhang]{\newline\indent School of Science
\newline\indent
Chongqing Jiaotong University
\newline\indent
Chongqing 400074, PR China
\newline\indent and
\newline\indent Chern Institute of Mathematics
\newline\indent
Nankai University
\newline\indent
Tianjin 300071, PR China}
\email{\href{mailto:zhangjianjun09@tsinghua.org.cn}{zhangjianjun09@tsinghua.org.cn}}

\address[J.M. do \'O]{\newline\indent Department of Mathematics
\newline\indent
Federal University of Para\'{\i}ba
\newline\indent
58051-900, Jo\~ao Pessoa-PB, Brazil}
\email{\href{mailto:jmbo@pq.cnpq.br}{jmbo@pq.cnpq.br}}

\address[M.\ Squassina]{\newline\indent Dipartimento di Informatica
\newline\indent
Universit\`a degli Studi di Verona,
\newline\indent
C\'a Vignal 2, Strada Le Grazie 15, I-37134 Verona, Italy}
\email{\href{mailto:marco.squassina@univr.it}{marco.squassina@univr.it}}

\thanks{Research partially supported by INCTmat/MCT/Brazil.\
J.M.\ do \'O was supported by CNPq, CAPES/Brazil,
Jianjun Zhang  was partially supported by CAPES/Brazil and CPSF (2013M530868)}


\begin{document}

\begin{abstract}
We consider a fractional Schr\"{o}dinger-Poisson system with a general nonlinearity in subcritical and critical case. The Ambrosetti-Rabinowitz condition
is not required. By using a perturbation approach, we prove the existence of positive solutions.
Moreover, we study the asymptotics of solutions for a vanishing parameter.
\end{abstract}
\maketitle

\s{Introduction and main result}

\noindent
We are concerned with the fractional nonlinear
Schr\"{o}dinger-Poisson system
\be\lab{q1} \left\{
\begin{array}{ll}
(-\Delta)^s u+\la\phi u=g(u)&\mbox{in}\ \R^3,\\
(-\Delta)^t \phi=\la u^2&\mbox{in}\ \R^3,
\end{array}
\right. \ee where $\la>0$, $(-\Delta)^\al$ is the fractional Laplacian operator for $\al=s,t\in (0,1)$. The fractional Schr\"{o}dinger equation was introduced by Laskin \cite{Lask3} and arose in fractional quantum mechanics in the study of particles on stochastic fields modeled by L\'{e}vy processes. The operator $(-\Delta)^\al$ can be seen as the infinitesimal generators of L\'{e}vy stable diffusion processes \cite{App}. If $\lambda=0$, the system \re{q1} reduces to the nonlinear fractional scalar field equation
\begin{equation}
\label{lask}
(-\Delta)^s u=g(u)\ \ \ \mbox{in $\R^3$.}
\end{equation}
This equation is related to the standing waves for fractional scalar field equation
\begin{equation}
\label{eq:SCH}
{\rm i}\phi_t-(-\Delta)^s\phi+g(\phi)=0 \ \ \ \text{in $\R^3$,}
\end{equation}
which is a physically relevant generalization of the classical NLS. For power type nonlinearities
the fractional Schr\"odinger equation was derived by Laskin \cite{Lask3} by replacing
the Brownian motion in the path integral approach with the so called L\'evy flights, see e.g.\  \cite{metkla1}. So, the equation
we want to study presents first of all as a perturbation of a physically meaningful equation.
Also, in \cite{FLS,FLS1} the author obtained deep results about uniqueness and non-degeneracy of ground states
for \eqref{lask} in case $g(u)=|u|^{p-2}u-u$ for subcritical $p$.
See also \cite{secsqu} where
the soliton dynamics for \eqref{eq:SCH} with an external potential was investigated.  In \cite{giammetta},
the author studies the evolution equation associated with the one dimensional system
\be\lab{q2} \left\{
\begin{array}{ll}
	-\Delta u+\la\phi u=g(u)&\mbox{in}\ \R,\\
	(-\Delta)^t \phi=\la u^2&\mbox{in}\ \R.
\end{array}
\right. \ee
In this case the diffusion is fractional only in the Poisson equation. Our system is more general
and contain this as a particular case.
If $\mathcal K_{\alpha}(x)=|x|^{\alpha-N},$
in \cite{boson} the following equation is studied
\begin{equation*}
	\label{bosonstar}
	\sqrt{-\Delta} u + u = (\mathcal{K}_2* |u|^2) u,
	\qquad
	u \in H^{1/2}(\mathbb{R}^3),\, \, u > 0,
\end{equation*}
and in \cite{ElSc} it is shown that  the dynamical evolution of boson stars is  described
by the nonlinear evolution equation
\begin{equation*}
	\label{bosonstarm}
	i\partial_t \psi= \sqrt{-\Delta + m^2} \psi - (\mathcal{K}_2* |\psi|^2) \psi \qquad (m\geq0)
\end{equation*}
for a field $\psi : [0, T ) \times \mathbb{R}^3 \to \mathbb{C}$
(see also \cite{FJL1}).
The square root of the Laplacian also appears in the semi-relativistic
Schr\"odinger-Poisson-Slater systems, \cite{JTN}. See also
the model studied in \cite{DSS}.
\noindent Observe that, taking formally $s=t=1$, then system \eqref{q1} reduces to the classical Schr\"odinger-Poisson system
\be\lab{sp} \left\{
\begin{array}{ll}
	-\Delta u+\la\phi u=g(u)&\mbox{in}\ \R^3,\\
	-\Delta \phi=\la u^2&\mbox{in}\ \R^3.
\end{array}
\right. \ee
It describes systems of identically charged particles interacting each other in the case where magnetic effects can be neglected \cite{Benci1}. In recent years, the Schr\"{o}dinger-Possion system \re{sp} has been widely studied by many researchers. Here we would like to cite some related results, for example, positive solutions \cite{CV}, ground state solutions \cite{AP}, semi-classical states \cite{T3}, sign-changing solutions \cite{Ianni1}. See also \cite{Am} and the references therein. In \cite{Azzo}, Azzollini, d'Avenia and Pomponio were concerned with \re{sp} under the Berestycki-Lions conditions $(H2)$-$(H4)$ with $s=1$. The authors proved that \re{sp} admits a positive radial solution if $\la>0$ small enough. For the critical case, we refer to \cite{Jian} and a recent work \cite{Zhang} of the authors of the present work.

\subsection{Main results} In the present paper, we are mainly concerned with the positive solutions of \re{q1}. First, we consider the subcritical case with the Berestycki-Lions conditions. Precisely, we assume the following hypotheses on $g$: \begin{itemize}
\item [$(H1)$] $g\in C^1(\mathbb{R},\mathbb{R})$;
\item [$(H2)$] $\displaystyle -\infty<\liminf_{\tau\rightarrow0}\frac{g(\tau)}{\tau}\le
\limsup_{\tau\rightarrow0}\frac{g(\tau)}{\tau}=-m<0$;
\item [$(H3)$] $\displaystyle
\limsup_{\tau\rightarrow\infty}\frac{g(\tau)}{\tau^{2^\ast_s-1}}\le0$, \, where $2^\ast_s=\frac{6}{3-2s}$;
\item [$(H4)$] there exists $\xi>0$ such that $G(\xi):=\int_0^\xi
g(\tau)\, \ud \tau>0$.
\end{itemize}
Our first result can read as
\bt\lab{Theorem 1} Suppose that $g$ satisfies $(H1)$-$(H4)$ and $2t+4s\ge3$.
\begin{itemize}
\item [(i)] There exists $\la_0>0$ such that, for every $\la\in (0,\la_0)$, system
\re{q1} admits a nontrivial positive radial solution $(u_\la,\phi_\la)$.
\item [(ii)] Along a subsequence, $(u_\la,\phi_\la)$ converges to $(u,0)$ in
$H^s(\R^3)\times\mathcal{D}^{t,2}(\R^3)$ as $\la\rg 0$, where $u$ is
a radial ground state solution of \re{lask}.
\end{itemize}
\et
\br\rm
The hypotheses $(H2)$-$(H4)$ are the so-called Berestycki-Lions conditions, which were introduced in \cite{Lions} to get the ground state of \re{lask} with $s=1$. Under $(H1)$-$(H4)$, X. Chang and Z.-Q. Wang \cite{Chang} proved the existence of ground state solutions to \re{lask} for $s\in(0,1)$. The hypothesis $(H1)$ is only used to get the better regularity of solutions to \re{lask}, which can guarantee the Pohoz\v{a}ev's identity. By the Pohoz\v{a}ev's identity, $(H4)$ is necessary.
\er
\br\rm
The hypothesis $2t+4s\ge3$ is just used to guarantee that the Poisson equation $(-\Delta)^t \phi=\la u^2$ make sense, due to $\mathcal{D}^{t,2}(\R^3)\hookrightarrow L^{2^\ast_t}(\R^3)$. For the details, see Section 2 below.
\er
\noindent In \cite{Zhang}, without the Ambrosetti-Rabinowtiz condition, the authors of the present work considered the existence and concentration of positive solutions to \re{q1} in the critical case for $s=t=1$. It is natural to wonder if similar results can hold for the critical fractional case. This is just our second goal of the present paper. In the critical case, we assume the following hypotheses on $g$:
\begin{itemize}
\item [$(H2)'$] $\displaystyle \lim_{\tau\rightarrow0}\frac{g(\tau)}{\tau}=-a<0$;
\item [$(H3)'$] $\displaystyle
\lim_{\tau\rightarrow\infty}\frac{g(\tau)}{\tau^{2^\ast_s-1}}=b>0$;
\item [$(H4)'$] there exists $\mu>0$ and $q<2^{\ast}_s$ such that $g(\tau)-b \tau^{2^{\ast}_s-1}+a\tau\ge \mu \tau^{q-1}$ for all $\tau>0$.
\end{itemize}
Our second result can read as follows.
\bt\lab{Theorem 2} Suppose that $g$ satisfies $(H1)$, $(H2)'$-$(H4)'$.
\begin{itemize}
\item [(i)] The limit problem \re{lask} admits a ground state solution if $\max\{2_s^\ast-2,2\}<q<2_s^\ast$.
\item [(ii)] Let $2t+4s\ge3$, then there exists $\la_0>0$ such that, for every $\la\in (0,\la_0)$, system
\re{q1} admits a nontrivial positive radial solution $(u_\la,\phi_\la)$ if $\max\{2_s^\ast-2,2\}<q<2_s^\ast$.
\item [(iii)] Along a subsequence, $(u_\la,\phi_\la)$ converges to $(u,0)$ in
$H^s(\R^3)\times\mathcal{D}^{t,2}(\R^3)$ as $\la\rg 0$, where $u$ is
a radial ground state solution of \re{lask}.
\end{itemize}
\et
\br\rm
In the case $s=1$, the hypotheses $(H2)'$-$(H4)'$ were introduced in \cite{Zhang-Zou} (see also \cite{Alves}) to obtain the ground state of the scalar field equation $-\DD u=g(u)$ in $\RN$. In \cite{Jihui}, X.\ Shang and J.\ Zhang considered the fractional problem \re{lask} in the critical case (see also \cite{Jihui1}). With the help of the monotonicity of $\tau\mapsto g(\tau)/\tau$, the ground state solutions were obtained by using the Nehari approach. To the best of our knowledge, there are few results in the literature about the ground state of the critical fractional problem \re{lask} with a general nonlinearity, particularly without the Ambrosetti-Rabinowtiz condition and the monotonicity of $g(\tau)/\tau$. Theorem \ref{Theorem 2} seems to be the first result
in this direction.
\er
\br\rm
Without loss generality, from now on, we assume that $a=b=\mu=1$.
\er
\noindent
In the rest of the paper, we use the perturbation approach to prove Theorem \ref{Theorem 1} and \ref{Theorem 2}. Similar argument also can be found in \cite{Zhang}.
\vskip0.1in
\noindent The paper is organized as follows. \newline
In Section 2 we introduce the functional framework and some preliminary results. \newline
In Section 3 we construct the min-max level. \newline
In Section 4, we use a perturbation argument to complete the proof of Theorem~\ref{Theorem 1}.\newline
In Section 5, we give the proof of Theorem~\ref{Theorem 2}.

\vskip0.1in
\noindent{\bf Notations.}
\begin{itemize}
\item [$\bullet$] $\|u\|_p:=\big(\int_{\R^3}|u|^p\, \ud x\big)^{1/p}$ for $p\in [1,\infty)$.
\item [$\bullet$] $2^\ast_\al:=\frac{6}{3-2\al}$ for any $\al\in(0,1)$.
\item [$\bullet$] $\hat{u}=\mathcal{F}(u)$ is the Fourier transform of $u$.
\end{itemize}
\vskip0.1in

\s{Preliminaries and functional setting}

\renewcommand{\theequation}{2.\arabic{equation}}

\subsection{Fractional order Sobolev spaces}
\noindent The fractional Laplacian $(-\Delta)^\al$ with $\al\in(0,1)$ of a function $\phi:\R^3\rightarrow\R$ is defined by
$$
\mathcal{F}((-\Delta)^\al\phi)(\xi)=|\xi|^{2\al}\mathcal{F}(\phi)(\xi),\, \, \xi\in\R^3,
$$
where $\mathcal{F}$ is the Fourier transform, i.e., $$\mathcal{F}(\phi)(\xi)=\frac{1}{(2\pi)^{3/2}}\int_{\R^3}\exp{(-2\pi i\xi\cdot x)}\phi(x)\, \ud x,$$ $i$ is the image unit. If $\phi$ is smooth enough, it can be computed by the following singular integral
$$
(-\Delta)^\al\phi(x)=c_\al\, \mbox{P.V.}\int_{\R^3}\frac{\phi(x)-\phi(y)}{|x-y|^{3+2\al}}\, \ud y,\, \ x\in\R^3,
$$
where $c_\al$ is a normalization constant and P.V. stands the principal value.

\noindent For any $\al\in(0,1)$, we consider the fractional order Sobolev space
$$
H^\al(\R^3)=\left\{u\in L^2(\R^3): \int_{\R^3}|\xi|^{2\al}|\hat{u}|^2\, \ud\xi<\iy\right\},
$$
endowed with the norm $$\|u\|_\al=\left(\int_{\R^3}(1+|\xi|^{2\al})|\hat{u}|^2\, \ud\xi\right)^{1/2},\,  u\in H^\al(\R^3),$$
and the inner product
$$
(u,v)_\al=\int_{\R^3}(1+|\xi|^{2\al})\hat{u}\bar{\hat{v}}\, \ud\xi,\,  u,v\in H^\al(\R^3).
$$
It is easy to know the inner products on $H^s(\R^3)$
\begin{itemize}
\item [] $u,v\mapsto \int_{\R^3}(1+|\xi|^{2\al})\hat{u}\bar{\hat{v}}\, \ud\xi$,
\vskip0.08in
\item [] $u,v\mapsto \int_{\R^3}\left(uv+(-\DD)^{\al/2}u(-\DD)^{\al/2}v\right)\, \ud x$,
\end{itemize}
are equivalent(see \cite{Jihui}).
The homogeneous Sobolev space $\mathcal{D}^{\al,2}(\R^3)$ is defined by
$$
\mathcal{D}^{\al,2}(\R^3)=\{u\in L^{2^\ast_\al}(\R^3): |\xi|^\al\hat{u}\in L^2(\R^3)\},
$$
which is the completion of $C_0^\iy(\R^3)$ under the norm
$$
\|u\|_{\mathcal{D}^{\al,2}}^2=\|(-\Delta)^{\al/2}u\|_2^2=\int_{\R^3}|\xi|^{2\al}|\hat{u}|^2\, \ud\xi,\,  u\in \mathcal{D}^{\al,2}(\R^3),
$$
and the inner product
$$
(u,v)_{\mathcal{D}^{\al,2}}=\int_{\R^3}(-\DD)^{\al/2}u(-\DD)^{\al/2}v\, \ud x,\,  u,v\in \mathcal{D}^{\al,2}(\R^3).
$$
For the further introduction on the Fractional order Sobolev space, we refer to \cite{DPV}. Let $$H_r^s(\R^3)=\{u\in H^3(\R^3): u(x)=u(|x|)\}.$$ Now, we introduce the following Sobolev embedding theorems.
\bl[see \cite{Lions1}]\lab{l1} For any $\al\in(0,1)$,
$H^\al(\R^3)$ is continuously embedded into $L^q(\R^3)$ for $q\in[2,2^\ast_\al]$ and compactly embedded into $L^q_{loc}(\R^3)$ for $q\in[1,2^\ast_\al)$. Moreover, $H_r^\al(\R^3)$ is compactly embedded into $L^q(\R^3)$ for $q\in(2,2^\ast_\al)$.
\el

\bl[see \cite{Cots,DPV}]\lab{l2} For any $\al\in(0,1)$,
$\mathcal{D}^{\al,2}(\R^3)$ is continuously embedded into $L^{2^\ast_\al}(\R^3)$, i.e., there exists $S_\al>0$ such that
$$
\left(\int_{\R^3}|u|^{2^\ast_\al}\, \ud x\right)^{2/2^\ast_\al}\le S_\al\int_{\R^3}|(-\Delta)^{\al/2}u|^2\, \ud x,\, \, u\in \mathcal{D}^{\al,2}(\R^3).
$$
\el
\subsection{The variational setting}

\noindent Now, we study the variational setting of \re{q1}. By Lemma \ref{l1}, $$H^s(\R^3)\hookrightarrow L^{12/(3+2t)}(\R^3)\, \, \mbox{if}\, \, 2t+4s\ge3.$$  Then, for $u\in H^s(\R^3)$, by Lemma \ref{l2} the linear operator $P:\mathcal{D}^{t,2}(\R^3)\rightarrow\R$ defined by
$$
P(v)=\int_{\R^3}u^2v\le\|u\|_{12/(3+2t)}^2\|v\|_{2^\ast_t}\le C\|u\|_s^2\|v\|_{\mathcal{D}^{t,2}},
$$
is well defined on $\mathcal{D}^{t,2}(\R^3)$ and is continuous. Thus, it follows from the Lax-Milgram theorem that there exists a unique $\phi_u^t\in
\mathcal{D}^{t,2}(\R^3)$ such that $(-\Delta)^t\phi_u^t=\la u^2$. Moreover, for $x\in\R^3$,
\be
\lab{t1} \phi_u^t(x):=\la c_t\int_{\R^3}\frac{u^2(y)}{|x-y|^{3-2t}}\, \ud y,
\ee
where we have set
$$
c_t=\frac{\G(\frac{3}{2}-2t)}{\pi^{\frac{3}{2}}2^{2t}\G(t)}.
$$
Formula \re{t1} is called the $t$-Riesz potential. Substituting \re{t1} into \re{q1}, we can rewrite
\re{q1} in the following equivalent form
\be\lab{q2} (-\DD)^s u+\la\phi_u^tu=g(u),\ \ \ u\in H^s(\R^3). \ee We
define the energy functional $\G_\la:H^s(\R^3)\to\R$ by
$$
\G_\la(u)=\frac{1}{2}\int_{\R^3}|(-\DD)^{s/2}u|^2\, \ud x+\frac{\la}{4}\int_{\R^3}\phi_u^tu^2\, \ud x-\int_{\R^3}G(u)\, \ud x,
$$
with $G(\tau)=\int_0^\tau g(\zeta)\, \ud \zeta$. Obviously, the critical points of $\G_\la$ are the weak solutions of \re{q2}.
\bd\noindent
\begin{itemize}
\item [$(1)$] We call $(u,\phi)\in H^s(\R^3)\times\mathcal{D}^{t,2}(\R^3)$ is a weak solution of \re{q1} if $u$ is a weak solution of \re{q2}.
\item [$(2)$] We call $u\in H^s(\R^3)$ is a weak solution of \re{q2} if
$$
\int_{\R^3}((-\DD)^{s/2}u(-\DD)^{s/2}v+\la\phi_u^tuv)\, \ud x=\int_{\R^3}g(v)v\, \ud x,\, \forall v\in H^s(\R^3).
$$
\end{itemize}
\ed
\noindent Setting
$$
T(u):=\frac{1}{4}\int_{\R^3}\phi_u^tu^2\, \ud x,
$$
then we summarize some properties of $\phi_u^t, T(u)$, which will be used later.

\bl
\lab{l3} If $t,s\in(0,1)$ and $2t+4s\ge3$, then for any $u\in H^s(\R^3)$, we have
\begin{itemize}
\item [(1)] $u\mapsto\phi_u^t: H^s(\R^3)\mapsto \mathcal{D}^{t,2}(\R^3)$ is continuous and maps bounded sets into bounded
sets.
\item [(2)] $\phi_u^t(x)\ge0, x\in\R^3$ and $T(u)\le
c\la\|u\|_s^4$ for some $c>0$.
\item [(3)] $T(u(\frac{\cdot}{\tau}))=\tau^{3+2t}T(u)$ for any $\tau>0$ and $u\in H^s(\R^3)$.
\item [(4)] If $u_n\rg u$ weakly in $H^s(\R^3)$, then $\phi_{u_n}\rg \phi_u$ weakly in
$\mathcal{D}^{t,2}(\R^3)$.
\item [(5)] If $u_n\rg u$ weakly in $H^s(\R^3)$, then
$T(u_n)=T(u)+T(u_n-u)+o(1)$.
\item [(6)] If $u$ is a radial function, so is $\phi_u^t$.
\end{itemize}
\el
\bp
The proof is similar as that in \cite{Ruiz}, so we omit the details here.
\ep
\vskip0.1in


\s{The subcritical case}
\renewcommand{\theequation}{3.\arabic{equation}}

\subsection{The modified problem}
\noindent It follows from Lemma \ref{l3} that $\G_\la$ is well defined on $H^s(\R^3)$ and is of class $C^1$. Since we are concerned with the positive solutions
of \re{q2}, similar as that in \cite{Lions}(see also \cite{Chang}), we modify our problem first. Without loss generality, we assume that $$0<\xi=\inf\{\tau\in(0,\iy): G(\tau)>0\},$$ where $\xi$ is given in $(H4)$.
Let $\tau_0=\inf\{\tau>\xi: g(\tau)=0\}\in[\xi,\iy]$, define $\ti{g}:\R\to\R$,
\begin{align*}
\ti{g}(\tau)=\left\{\begin{array}{ccc}
  g(\tau)&\mbox{if}\, \tau\in[0,\tau_0] \\
  0&\mbox{if}\, \tau\ge\tau_0,
\end{array}\right.
\end{align*}
and $\ti{g}(\tau)=0$ for $\tau\le 0$. If $u\in H^s(\R^3)$ is a solution of \re{q2} where $g$ is replaced by $\ti{g}$, then by the maximum principle \cite{Sire} we get that $u$ is positive and $u(x)\le\tau_0$ for any $x\in\R^3$, i.e., $u$ is a solution of the original problem \re{q2} with $g$. Thus, from now on, we can replace $g$ by $\ti{g}$, but still use the same notation $g$.
In addition, for $\tau>0$, let
$$
g_1(\tau)=\max\{g(\tau)+m\tau,0\},\, \, g_2(\tau)=g_1(\tau)-g(\tau),
$$
then $g_2(\tau)\ge m\tau$ for $\tau\ge 0$,
\be\lab{s1}
\lim_{\tau\rg0}\frac{g_1(\tau)}{\tau}=0,\, \, \lim_{\tau\rg+\iy}\frac{g_1(\tau)}{\tau^{2^\ast_s-1}}=0,
\ee
and for any $\e>0$, there exists $C_\e>0$ such that
\be\lab{s2}
g_1(\tau)\le\e g_2(\tau)+C_\e\tau^{2^\ast_s-1},\, \tau\ge0.
\ee
Let $G_i(u)=\int_0^ug_i(\tau)\, \ud\tau, i=1,2$, then by \re{s1} and \re{s2} for any $\e>0$ there exists $C_\e>0$ such that
\be\lab{s3}
G_1(\tau)\le\e G_2(\tau)+C_\e|\tau|^{2^\ast_s},\, \tau\in\R.
\ee
\subsection{The limit problem} In the following, we will find the solutions of \re{q2} by seeking the critical points of $\G_\la$. If $\la=0$, problem \re{q2} becomes
\be
\lab{q3} (-\DD)^su=g(u),\ \ \ u\in H^s(\R^3),
\ee
which is referred as the limit problem of \re{q2}. We define an energy functional for the limiting problems \eqref{q3} by
$$
L(u)=\frac{1}{2}\int_{\R^3}|(-\DD)^{s/2}u|^2\, \ud x-\int_{\R^3}G(u)\, \, \ud x,\ \ u\in H^s(\R^3).
$$
In \cite{Chang}, X. Chang and Z.-Q. Wang proved that, with the same assumptions on $g$ in Theorem \ref{Theorem 1}, there exists a positive ground state solution $U\in H_r^3(\R^3)$ of \eqref{q3}. Moreover, each such solution $U$
of \eqref{q3} satisfies the Pohoz\v{a}ev identity
\be\label{eq}
\frac{3-2s}{2}\int_{\R^3}|(-\DD)^{s/2}U|^2\, \ud x=3\int_{\R^3}G(U)\, \ud x.
\ee
Let $S$ be the set of positive radial ground state solutions $U$ of \eqref{q3}, then $S\not=\phi$ and we have the following compactness result, which plays a crucial role in the proof of Theorem \ref{Theorem 1}.
\bo\lab{prop2.1}\
Under the assumptions in Theorem \ref{Theorem 1}, $S$ is compact in $H_r^s(\R^3)$.
\eo
\noindent As shown in \cite{Cho}, for general $s\in(0,1)$, we do not have a similar radial lemma in $H_r^s(\R^3)$. So the Strauss's compactness lemma(see \cite{Lions}) is not applicable here. Before we prove Proposition \ref{prop2.1}, we start with the following compactness lemma, which is a special case of \cite[Lemma 2.4.]{Chang}
\bl[see \cite{Chang}]\lab{l4}
Assume $Q\in C(\R,\R)$ satisfy
$$
\lim_{\tau\rg0}\frac{Q(\tau)}{\tau^2}=\lim_{|\tau|\rg\iy}\frac{Q(\tau)}{|\tau|^{2^\ast_s}}=0,
$$
and there exists a bounded sequence $\{u_n\}_{n=1}^\iy\subset H_r^s(\R^3)$ for some $v\in L^1(\R^3)$ with
$$
\lim_{n\rg\iy}Q(u_n(x))=v(x),\, a.e.\, x\in\R^3.
$$
Then, up to a subsequence, we have $Q(u_n)\rightarrow v$ strongly in $L^1(\R^3)$ as $n\rg\iy$.
\el
\noindent{\it Proof of Proposition \ref{prop2.1}} Let $\{u_n\}_{n=1}^\iy\subset S$ and denote by $E$ the least energy of \re{q3}, then for any $n$, $u_n$ satisfies $L(u_n)=E$ and the Pohoz\v{a}ev identity \re{eq}, which implies that
$$
E=\frac{s}{3}\int_{\R^3}|(-\DD)^{s/2}u_n|^2\, \ud x\, \, \mbox{and}\, \, \int_{\R^3}G(u_n)\, \ud x=\frac{3-2s}{2s}E.
$$
Obviously, $\{\|(-\DD)^{s/2}u_n\|_2\}$ is bounded. It follows from Lemma \ref{l2} that $\{\|u_n\|_{2^\ast_s}\}$ is bounded. By \re{s3}, as we can see in \cite{Lions}, $\{\|u_n\|_2\}$ is bounded, which yields that $\{u_n\}$ is bounded in $H_r^s(\R^3)$. Without loss generality, we can assume that there exists $u_0\in H_r^s(\R^3)$ such that $u_n\rg u_0$ weakly in $H_r^s(\R^3)$, strongly in $L^q(\R^3)$ for $q\in(2,2^\ast_s)$ and $u_n(x)\rg u_0(x)$ a.e. $x\in\R^3$.

\noindent In the following, we adopt some ideas in \cite{Lions} to prove that $u_n\rg u_0$ strongly in $H_r^s(\R^3)$. For $u\in H^s(\R^3)$, let
$$
J(u)=\frac{s}{3}\int_{\R^3}|(-\DD)^{s/2}u|^2\, \ud x\, \, \mbox{and}\, \, V(u)=\int_{\R^3}G(u)\, \ud x,
$$
then, we know $u_n$ is a minimizer of the following constrained minimizing problem
$$
\inf\left\{J(u): u\in H_r^s(\R^3), V(u)=\frac{3-2s}{2s}E\right\}.
$$
By \re{s1} and Lemma \ref{l4}, we get that
$$
\lim_{n\rg\iy}\int_{\R^3}G_1(u_n)=\int_{\R^3}G_1(u_0).
$$
Then by the Fatou's Lemma, $V(u_0)\ge\frac{3-2s}{2s}E$, which implies that $u_0\not\equiv0$. Meanwhile, it is easy to know that $J(u_0)\le E$. Similar as that in \cite{Lions}, we know that $u_0$ satisfies $J(u_0)=E$ and $V(u_0)=\frac{3-2s}{2s}E$. It yields that
$$
\lim_{n\rg\iy}\int_{\R^3}G_2(u_n)=\int_{\R^3}G_2(u_0).
$$
By the Fatou's Lemma, we know $\|u_n\|_2\rg\|u_0\|_2$ as $n\rg\iy$. Thus, $u_n\rg u_0$ strongly in $H_r^s(\R^3)$. The proof is completed.
\qed

\subsection{The minimax level}

\noindent Take $U\in S$, let $$U_\tau(x)=U(\frac{x}{\tau}), \tau>0,$$ then by the definition of $\hat{U}=\mathcal{F}(U)$, we know
$\hat{U}(\frac{\cdot}{\tau})=\tau^3\hat{U}(t\cdot)$. Then
$$
\int_{\R^3}|(-\DD)^{s/2}U_\tau|^2\, \ud x=\int_{\R^3}|\xi|^{2s}|\hat{U}(\frac{\xi}{\tau})|^2=\tau^{3-2s}\int_{\R^3}|(-\DD)^{s/2}U|^2\, \ud x.
$$
By the Pohoz\v{a}ev's identity,
$$
L(U_\tau)=\Big(\frac{\tau^{3-2s}}{2}-\frac{3-2s}{6}\tau^3\Big)\int_{\R^3}|(-\DD)^{s/2}U|^2.
$$
Thus, there exists $\tau_0>1$ such that $L(U_{\tau})<-2$ for $\tau\ge \tau_0$. Set
$$
D_\la\equiv\max_{\tau\in[0,\tau_0]}\G_\la(U_\tau).
$$
By virtue of Lemma \ref{l3}, $\G_\la(U_\tau)=L(U_{\tau})+O(\la)$. Note that $\max_{\tau\in[0,\tau_0]}L(U_\tau)=E$, we get that $D_\la\rg E$, as $\la\rg 0^+$.

\noindent
Moreover, similar to \cite{Zhang}, we can prove the following lemma, which is crucial to define the uniformly bounded set of the mountain pathes(see below).
\bl\lab{l4.1} There exist $\la_1>0$ and  $\mathcal{C}_0>0$, such that for
any $0<\la<\la_1$ there hold
\begin{equation*}
\G_\la(U_{\tau_0})<-2,\qquad
\|U_\tau\|_s\le \mathcal{C}_0, \,\,\, \forall \tau\in (0,\tau_0],\qquad
\|u\|_s\le \mathcal{C}_0, \,\,\, \forall u\in S.
\end{equation*}
\el
\noindent Now, for any $\la\in (0,\la_1)$, we define a min-max value $C_\la$:
$$C_\la=\inf_{\g\in \Upsilon_\la}\max_{\tau\in [0,\tau_0]}\G_\la(\g(\tau)),$$
where
\begin{align*}
\Upsilon_\la=\big\{\g\in C([0,\tau_0],H_r^s(\R^3)): \g(0)=0,
\g(\tau_0)=U_{\tau_0},\|\g(\tau)\|_s\le \mathcal{C}_0+1,\tau\in[0,\tau_0]\big\}.
\end{align*}
Obviously, for $\tau>0$, $$\|U_\tau\|_s^2=\tau^{3-2s}\|(-\DD)^{s/2}U\|_2^2+\tau^3\|U\|_2^2.$$ Then we can define $U_0\equiv0$, so $U_\tau\in \Upsilon_\la$. Moreover,
$$\limsup_{\la\rg0^+}C_\la\le\lim_{\la\rg0^+}D_\la=E.$$

\bo\lab{bo4} $\lim\limits_{\la\rg 0^+}C_\la=E.$ \eo

\bp It suffices to prove that
$$
\liminf\limits_{\la\rg 0^+}C_\la\ge E.
$$
Now, we give the following mountain pass value
$$b=\inf_{\g\in \Upsilon}\max_{\tau\in [0,1]}L(\g(\tau)),$$
where $\Upsilon=\big\{\g\in C([0,1],H_r^s(\R^3)): \g(0)=0,\g(1)<0\big\}.$ It follows from \cite[Lemma 3.2]{Chang} that $L$ satisfies the mountain pass geometry. As we can see in \cite{Jean}, $b$ agrees with the least energy level of \re{q3}, i.e., $b=E$. Note that $\phi_u^t(x)\ge 0, x\in\R^3$, then $\ti{\g}(\cdot)=\g(\tau_0\cdot)\in \Upsilon$ for any $\g\in \Upsilon_\la$. It follows that $C_\la\ge
b$, concluding the proof. \ep

\subsection{Proof of Theorem \ref{Theorem 1}}

\noindent
Now for $\al,d>0$, define
$$
\G_\la^\al:=\{u\in H_r^s(\R^3): \G_\la(u)\le\al\}
$$
and
$$
S^d=\left\{u\in H_r^s(\R^3):  \inf_{v\in
S}\|u-v\|_s\le d\right\}.
$$
In the following, we will find a solution $u\in S^d$ of problem \re{q2} for
sufficiently small $\la>0$ and some $0<d<1$. The following proposition is crucial to obtain a suitable
(PS)-sequence for $\G_\la$ and plays a key role in our proof.

\bo\lab{bo5} Let $\{\la_i\}_{i=1}^\infty$ be such that $\lim_{i\rg
\infty}\la_i=0$ and $\{u_{\la_i}\}\subset S^d$
with
$$
\lim_{i\rg\infty}\G_{\la_i}(u_{\la_i})\le E\ \mbox{and}\
\lim_{i\rg\infty}\G_{\la_i}^{'}(u_{\la_i})=0.
$$
Then for $d$ small enough, there is $u_0\in S$, up to
a subsequence, such that $u_{\la_i}\rg u_0$ in $H_r^s(\R^3)$.
\eo
\bp
For convenience, we write $\la$ for $\la_i$. Since $u_\la\in
S^d$ and $S$ is compact, we know $\{u_\la\}$ is bounded in $H_r^s(\R^3)$. Then by Lemma \ref{l3} we see that
$$
\lim_{i\rg\infty}L(u_{\la})\le E\ \mbox{and}\
\lim_{i\rg\infty}L^{'}(u_{\la})=0.
$$
It follows from \cite[Lemma 3.3]{Chang} that, there is $u_0\in H_r^s(\R^3)$, up to
a subsequence, such that $u_\la\rg u_0$ strongly in $H_r^s(\R^3)$. Obviously, $0\not\in S^d$ for $d$ small. This implies that $u_0\not\equiv0$, $L(u_0)\le E$ and $L'(u_0)=0$. Thus, $L(u_0)=E$, i.e., $u_0\in S$. The proof is completed. \ep

\noindent
By Proposition \ref{bo5}, for small $d\in(0,1)$ such that there exist $\omega>0, \la_0>0$ such that 
\be\lab{ge1}\hbox{$\|\G_\la'(u)\|_s\ge\omega$ for
$u\in\G_\la^{D_\la}\bigcap(S^d\setminus S^{\frac{d}{2}})$
and $\la\in(0,\la_0)$.}
\ee
Similar to \cite{Zhang}, we have
\bo\lab{bo6} There exists $\al>0$ such that for small $\la>0$,
$$
\G_\la(\g(\tau))\ge C_\la-\al\ \ \mbox{implies that}\ \ \g(\tau)\in
S^{\frac{d}{2}},
$$
where $\g(\tau)=U(\frac{\cdot}{\tau}), \tau\in (0,\tau_0]$. \eo

\bp From Lemma \ref{l3} and the Pohoz\v{a}ev's identity,

$$\G_\la(\g(\tau))=\Big(\frac{\tau^{3-2s}}{2}-\frac{3-2s}{6}\tau^3\Big)\int_{\R^3}|(-\DD)^{s/2}U|^2+\la\tau^{3+2t}T(U).$$
Then $$\lim_{\la\rg0^+}\max_{\tau\in[0,\tau_0]}\G_\la(\g(\tau))=\max_{\tau\in[0,\tau_0]}\Big(\frac{\tau^{3-2s}}{2}-\frac{3-2s}{6}\tau^3\Big)\int_{\R^3}|(-\DD)^{s/2}U|^2=E.$$
The conclusion follows. \ep

\noindent
Similarly as that in \cite{Zhang}, thanks to \re{ge1} and Proposition \ref{bo6}, we can prove the following proposition, which assures the existence of a bounded
Palais-Smale sequence for $\G_\la$.

\bo\lab{bo7} For $\la>0$ small enough, there exists
$\{u_n\}_n\subset \G_\la^{D_\la}\cap S^d$ such that
$\G_\la'(u_n)\rg 0$ as $n\rg\iy$.
\eo

\vskip2pt
\noindent

\noindent{\bf Proof of Theorem \ref{Theorem 1} concluded.} It follows from Proposition \ref{bo7} that there exists
$\la_0>0$ such that for $\la\in(0,\la_0)$, there exists $\{u_n\}\in
\G_\la^{D_\la}\cap S^d$ with $\G_\la'(u_n)\rg 0$ as
$n\rg\iy$. Noting that $S$ is compact in $H_r^s(\R^3)$, we get that $\{u_n\}$ is bounded in $H_r^s(\R^3)$. Assume that $u_n\rg u_\la$ weakly in $H_r^s(\R^3)$, then
$\G_\la'(u_\la)=0$. It follows from the compactness of $S$ that $u_\la\in S^d$ and $\|u_n-u_\la\|_s\le 3d$ for $n$ large. So $u_\la\not\equiv0$ for small $d>0$.
By Lemma \ref{l3},
$$
\G_\la(u_n)=\G_\la(u_\la)+\G_\la(u_n-u_\la)+o(1).
$$
Noting that $G_2(\tau)\ge \frac{m}{2}\tau^2$ for any $\tau\in\R$, it follows from \re{s3} that for some $C>0$,
\begin{align*}
\G_\la(u_n-u_\la)\ge&\frac{1}{2}\int_{\R^3}(|(-\DD)^{s/2}(u_n-u_\la)|^2+\frac{m}{4}|u_n-u_\la|^2)\, \ud x\\
&-C\int_{\R^3}|u_n-u_\la|^{2^\ast_s}\, \ud x.
\end{align*}
Then by Lemma \ref{l2}, for small $d>0$, it is easy to verify that $\G_\la(u_n-u_\la)\ge0$ for large $n$. So $u_\la\in
\G_\la^{D_\la}\cap S^d$ with $\G_\la'(u_\la)=0$. Thus $u_\la$ is a nontrivial solution of \re{q2}. Finally, by Proposition \ref{bo5} we can get the asymptotic behavior of $u_\la$ as
$\la\rg 0^+$. The proof is completed. \qed
\vskip2pt

\s{The critical case}
\renewcommand{\theequation}{4.\arabic{equation}}
\noindent In this section, we consider the Schr\"{o}dinger-Possion system \re{q1} in the critical case. First, we establish the existence of ground state solutions to the fractional scalar field equation \re{lask} with a general critical nonlinear term. Then by the perturbation argument, we seek the solutions of \re{q1} in some neighborhood of the ground states to \re{lask}.
\vskip0.1in
\subsection{The limit problem} In this subsection, we use the constraint variational approach to seek the ground state solutions of \re{lask}. Similar argument also can be found in \cite{Lions,Zhang-Zou,Dip}. Let
$$T(u)=\frac{1}{2}\int_{\R^3}|(-\DD)^{s/2}u|^2\, \ud x, \,\,\,\quad  V(u)=\int_{\R^3}G(u)\, \ud x.$$
We recall that $U$ is said to be a ground state solution of \re{lask} if and only if $I(U)=m_0$, where $m_0:=\inf\{I(u):u\in H^s(\R^3)\setminus\{0\}\mbox{ is a solution of}\  \re{lask}\}$ and $$I(u)=T(u)-V(u).$$ The existence of ground state is reduced to look at the constraint minimization problem
\be\lab{q2n}
M:=\inf\{T(u): V(u)=1,u\in H^s(\R^3)\}
\ee
and eventually to remove the Lagrange multiplier by some appropriate scaling. Now, we state the main result in this subsection.
\bt\lab{gr}
Let $s\in(0,1)$ and assume that $(H2)'$-$(H4)'$ and
\begin{itemize}
\item [$(H0)$] $g\in C(\R,\R)$ and $g$ is odd, i.e., $g(-\tau)=-g(\tau)$ for $\tau\in\R$.
\end{itemize}
Then \re{lask} admits a positive ground state solution.
\et
\br\rm
Since we are concerned the positive solution of \re{lask}, $(H0)$ can be replaced by $(H0)':$ $g\in C(\R^+,\R)$. Moreover, similar to Theorem \ref{gr}, a similar result in $\RN(N>2s)$ also can be obtained.
\er
\noindent{\it Proof of Theorem \ref{gr}.} The proof follows the lines of that in \cite{Zhang-Zou}. For the completeness, we give the details here. \\
\noindent {\bf Step 1.} Let $M$ be given by \re{q2n} and $S_s$ be the Sobolev best constant in Lemma \ref{l2} for $s\in(0,1)$, then we claim that $$0<M<\frac{1}{2}\big(2_s^{\ast}\big)^{\frac{3-2s}{3}}S_s.$$
First, we prove that $\{u\in H^s(\R^3):V(u)=1\}\not=\phi$. By \cite{Best,Ser}, $S_s$ can be achieved by
$$
U_{\varepsilon}(x)=\kappa\e^{-\frac{3-2s}{2}}{\left(\mu^2+\left|\frac{x}{\e S_s^{\frac{1}{2s}}}\right|^2\right)^{-\frac{3-2s}{2}}},
$$
for any $\e>0$, where $\kappa\in\R, \mu>0$ are fixed constants. Let $\varphi\in C_0^{\infty}(\R^3)$ is a cut-off function with support $B_2$ such that $\varphi \equiv 1$ on $B_1$ and $0\le \varphi \le 1$ on $B_2$, where $B_r:=\{x\in\R^3: |x|<r\}$. Let $\psi_\varepsilon(x)=\varphi(x)U_{\varepsilon}(x)$, it follows from \cite{Ser} that
\be\lab{e1}
\int_{\R^3}|\psi_{\varepsilon}|^{2_s^{\ast}}=S_s^{3/(2s)}+O(\varepsilon^3), \quad \int_{\R^3}|(-\DD)^{s/2}\psi_{\varepsilon}|^2=S_s^{3/(2s)}+O(\varepsilon^{3-2s}).
\ee
Let
$$
v_\varepsilon=\frac{\psi_\varepsilon}{\|\psi_\varepsilon\|_{2_s^{\ast}}},
$$ then $\|(-\DD)^{s/2} v_\varepsilon\|_2^2\le S_s+O(\varepsilon^{3-2s})$.
Let
$$
\Gamma_\e:= \frac{1}{q}\|v_{\e}\|_q^q-\frac{1}{2}\|v_{\e}\|_2^2,
$$
then by $(H4)'$ we have
$V(v_\varepsilon)\ge \frac{1}{2_s^{\ast}}+\Gamma_\varepsilon$. In the following, we will show that
\be\lab{e4}
\lim_{\varepsilon\rightarrow 0}\frac{\Gamma_\varepsilon}{\varepsilon^{3-2s}}=+\infty.
\ee
By $\max\{2_s^\ast-2,2\}<q<2_s^\ast$, we know that $(3-2s)q>3$. Then it is easy to know there exist $C_1(s), C_2(s)>0$ such that
\begin{align*}
& \|v_\varepsilon\|_q^q\ge\frac{1}{\|\psi_\varepsilon\|_{2_s^{\ast}}^q}\int_{B_1}|U_\varepsilon|^q
\ge  C_1(s)\varepsilon^{3-\frac{3-2s}{2}q}\int_0^{\frac{1}{\e S_s^{1/(2s)}}}\frac{r^2}{(\mu^2+r^2)^{\frac{3-2s}{2}q}}dr\\
&\ \ \ \ \ \ \ \ =O(\varepsilon^{3-\frac{(3-2s)q}{2}}), \\
& \|v_\varepsilon\|_2^2 \le \frac{1}{\|\psi_\varepsilon\|_{2_s^{\ast}}^2}\int_{B_2}|U_\varepsilon|^2
\le C_2(s)\varepsilon^{2s}\int_0^{\frac{2}{\e S_s^{1/(2s)}}}\frac{r^2}{(\mu^2+r^2)^{3-2s}}dr\\
&\ \ \ \ \ \ \ \ =\left\{
                                 \begin{array}{ll}
                                   O(\e^{2s}), & \hbox{if } s<\frac{3}{4}; \\
                                   O(\e^{2s}\ln{\frac{1}{\e}}), & \hbox{if } s=\frac{3}{4}; \\
                                   O(\e^{3-2s}), & \hbox{if } s>\frac{3}{4}.
                                 \end{array}
                               \right.
\end{align*}
Then, we obtain
$$
\G_\e\ge O(\varepsilon^{3-\frac{(3-2s)q}{2}}),\quad \hbox{if } s\in(0,1).
$$
Noting that $\max\{2_s^\ast-2,2\}<q<2_s^\ast$, it is easy to verify that \re{e4} is true.
Thus, it follows that $V(v_\e)>0$ for small $\e>0$. By a scaling, we get that $\{u\in H^s(\R^3):V(u)=1\}\not=\phi$.
\vskip0.1in
\noindent Next, obviously, $M\in(0,+\iy)$. For small $\e>0$, $V(v_\e)>0$, then
\begin{eqnarray*}
M \le\frac{T(v_{\e})}{(V(v_{\e}))^{\frac{2}{2_s^{\ast}}}}\le\frac{1}{2}\frac{\|(-\DD)^{s/2} v_{\e}\|_2^2}{\left(\frac{1}{2_s^{\ast}}+\G_\e\right)^{\frac{2}{2_s^{\ast}}}}\le \frac{1}{2}(2_s^{\ast})^{\frac{2}{2_s^{\ast}}}S_s
\frac{1+O(\e^{N-2s})}{\left(1+2_s^{\ast} \G_\e\right)^{\frac{2}{2_s^{\ast}}}} .
\end{eqnarray*}
If $p\ge 1$, then $(1+t)^p\le 1+p(1+t)^{1+p}t$, for all $t\ge-1$. It follows from $(\ref{e4})$ that
$$
\left(1+O(\varepsilon^{N-2s})\right)^{\frac{2_s^{\ast}}{2}}-1\le \frac{2_s^{\ast}}{2}(1+O(\varepsilon^{N-2s}))^{1+\frac{2_s^{\ast}}{2}}O(\varepsilon^{N-2s})<2_s^{\ast}\Gamma_\varepsilon,
$$
for small $\varepsilon>0$, which yields $1+O(\varepsilon^{N-2s})<(1+2_s^{\ast}\Gamma_\varepsilon)^{2/{2_s^{\ast}}}$. Then $M<\frac{1}{2}\big(2_s^{\ast}\big)^{\frac{3-2s}{3}}S_s$.
\vskip0.1in
\noindent {\bf Step 2.} $M$ can be achieved. Noting that $g$ is odd and the Fractional Polya-Szeg\"{o} inequality \cite{Park}, without loss of generality, we can assume that there exists a positive minimizing sequence $\{u_n\}\subset H^s_r(\R^3)$ such that $V(u_n)=1$ and $T(u_n)\rightarrow M$ as $n\rightarrow \infty$. By Lemma \ref{l2}, it is easy to know that $\{u_n\}$ is bounded in $H_r^s(\R^3)$. By Lemma \ref{l1}, we can assume that $u_n\rg u_0$ weakly in $H^s(\R^3)$, strongly in $L^q(\R^3)$ and a.e. in $\R^3$. Set $v_n=u_n-u_0$, then $T(u_n)=T(v_n)+T(u_0)+o(1),$ and
$$
\|u_n\|_{2_s^{\ast}}^{2_s^{\ast}}=\|v_n\|_{2_s^{\ast}}^{2_s^{\ast}}+\|u_0\|_{2_s^{\ast}}^{2_s^{\ast}}+o(1),\, \,
\|u_n\|_2^2=\|v_n\|_2^2+\|u_0\|_2^2+o(1),
$$
where $o(1)\rightarrow 0$ as $n\rightarrow \infty$. Let $f(s)=g(s)-s^{2_s^{\ast}-1}+s$, then it follows from Lemma \ref{l4} that
$$
\int_{\R^3}F(u_n)=\int_{\R^3}F(u_0)+\int_{\R^3}F(v_n)+o(1).
$$
So, $V(u_n)=V(v_n)+V(u_0)+o(1)$.
\vskip0.1in
\noindent Next, we prove $u_0$ is the minimizer for $M$. Set $S_n=T(v_n), S_0=T(u_0), V(v_n)=\lambda_n, V(u_0)=\lambda_0$, we have $\lambda_n=1-\lambda_0+o(1)$ and $S_n=M-S_0+o(1)$. Under a scale change, we get that
\be\lab{q4n}
T(u)\ge M(V(u))^{\frac{3-2s}{3}},
\ee
for all $u\in H^s(\R^3)$ and $V(u)\geq 0$. By \re{q4n} we have $\lambda_0\in [0,1]$. If $\lambda_0\in (0,1)$, then by \re{q4n} again,
\begin{align*}
M&=\lim_{n\rightarrow \infty}(S_0+S_n)\nonumber\\
&\ge\lim_{n\rightarrow \infty}M\left((\lambda_0)^{\frac{3-2s}{3}}+(\lambda_n)^{\frac{3-2s}{3}}\right)\nonumber\\
&= M\left((\lambda_0)^{\frac{3-2s}{3}}+(1-\lambda_0)^{\frac{3-2s}{3}}\right)\nonumber\\
&> M(\lambda_0+1-\lambda_0)=M,
\end{align*}
which is a contradiction. On the other hand, if $\lambda_0=0$, then $S_0=0$, which implies $u_0=0$.
Then
$$
\limsup_{n\rightarrow\infty}\|v_n\|_{2_s^{\ast}}^2\ge (2_s^{\ast})^{\frac{3-2s}{3}}
$$
and
$$
M=\frac{1}{2}\lim_{n\rightarrow \infty}\|(-\DD)^{s/2}v_n\|_2^2\ge \frac{1}{2}(2_s^{\ast})^{\frac{3-2s}{3}}\liminf_{n\rightarrow\infty}\frac{\|(-\DD)^{s/2}v_n\|_2^2}{\|v_n\|_{2_s^{\ast}}^2}\ge\frac{1}{2}(2_s^{\ast})^{\frac{3-2s}{3}}S_s.
$$
which is again a contradiction. Then we conclude  $\lambda_0=1$, i.e., $M$ is achieved by $u_0$.

\vskip0.1in \noindent Finally, let $U(\cdot)=u_0(\cdot/\si_0)$, where $\si_0=\left(\frac{3-2s}{3}M\right)^{1/2}$, then $U$ is a ground state solution of \re{lask}. The proof is completed.
\qed

\br\rm \noindent Furthermore, similar as that in \cite{Chang}, if we assume that $g\in C^1(\R,\R)$ additionally, $U$ satisfies the Pohoz\v{a}ev identity
\be\label{eq}
\frac{3-2s}{2}\int_{\R^3}|(-\DD)^{s/2}U|^2\, \ud x=3\int_{\R^3}G(U)\, \ud x.
\ee
\noindent Similar as that in \cite{Jean,Zhang-Zou}, $U$ is also a mountain pass solution.
\er
\noindent Let $S_1$ be the set of positive radial ground state solutions $U$ of \eqref{lask}, then similar to Step 2 in the proof of Theorem \ref{gr}, we have the following compactness result.
\bo\lab{prop4.1}\
Under the assumptions in Theorem \ref{gr}, $S_1$ is compact in $H_r^s(\R^3)$.
\eo
\vskip0.1in
\subsection{Proof of Theorem \ref{Theorem 2}}

\noindent In the following, we are ready to prove Theorem \ref{Theorem 2}. Similar to Section 3, take $U\in S_1$, let $$U_\tau(x)=U(\frac{x}{\tau}), \tau>0,$$ then there exists $\tau_1>1$ such that $I(U_{\tau})<-2$ for $\tau\ge \tau_1$. Set
$$
D_\la^1\equiv\max_{\tau\in[0,\tau_1]}\G_\la(U_\tau),
$$
there exist $\la_2>0$ and $\mathcal{C}_1>0$ such that for
any $0<\la<\la_2$,
\begin{align*}
\phi\not=\Upsilon_\la=\big\{\g\in C([0,\tau_1],H_r^s(\R^3)): \g(0)=0,
\g(\tau_1)=U_{\tau_1},\|\g(\tau)\|_s\le \mathcal{C}_1+1,\tau\in[0,\tau_1]\big\}.
\end{align*}
Then for any $\la\in (0,\la_1)$, we define a min-max value $C_\la^1$:
$$C_\la^1=\inf_{\g\in \Upsilon_\la}\max_{\tau\in [0,\tau_1]}\G_\la(\g(\tau)).$$
Similar to Section 3, we have
\bo\lab{bo4n} $\lim\limits_{\la\rg 0^+}C_\la^1=\lim\limits_{\la\rg 0^+}D_\la^1=m,$ where $m$ is the least energy of \re{lask}.  \eo

\noindent
Now for $\al,d>0$, define
$$
\G_\la^\al:=\{u\in H_r^s(\R^3): \G_\la(u)\le\al\}
$$
and
$$
S_1^d=\left\{u\in H_r^s(\R^3):  \inf_{v\in
S_1}\|u-v\|_s\le d\right\}.
$$
Similar as that in Section 3, for small $\la>0$ and some $0<d<1$, we will find a solution $u\in S_1^d$ of \re{q2} in the critical case. Similar to \cite{Zhang}, we can get the following compactness result, which can yield the gradient estimate of $\G_\la$.

\bo\lab{bo5n} Let $\{\la_i\}_{i=1}^\infty$ be such that $\lim_{i\rg
\infty}\la_i=0$ and $\{u_{\la_i}\}\subset S_1^d$
with
$$
\lim_{i\rg\infty}\G_{\la_i}(u_{\la_i})\le m\ \mbox{and}\
\lim_{i\rg\infty}\G_{\la_i}^{'}(u_{\la_i})=0.
$$
Then for $d$ small enough, there is $u_1\in S_1$, up to
a subsequence, such that $u_{\la_i}\rg u_1$ in $H_r^s(\R^3)$.
\eo
\bp
For convenience, we write $\la$ for $\la_i$. Since $u_\la\in
S_1^d$ and $S_1$ is compact, we know $\{u_\la\}$ is bounded in $H_r^s(\R^3)$. Moreover, up to a subsequence, there exists $u_1\in
S_1^d$ such that $u_\la\rg u_1$ weakly in $H^s(\R^3)$, a.e. in
$\R^3$ and $\|u_\la-u_1\|_s\le3d$ for $i$ large. Then by Lemma \ref{l3} we see that
$$
\lim_{i\rg\infty}I(u_{\la})\le m\ \mbox{and}\
\lim_{i\rg\infty}I^{'}(u_{\la})=0.
$$
Then $I'(u_1)=0$. Obviously, $u_0\not\equiv 0$ if $d$ small. So $I(u_1)\ge m$. Meanwhile, thanks to
Lemma \ref{l4}, $I(u_\la)=I(u_1)+I(u_\la-u_1)+o(1)$, then we
have $$I(u_\la-u_1)=\frac{1}{2}\|u_\la-u_1\|_s^2-\frac{1}{2_s^\ast}\|u_\la-u_1\|_{2_s^\ast}^{2_s^\ast}+o(1)\le o(1).$$
Then by Lemma \ref{l2}, for $d$ small enough, $u_\la\rg u_1$ strongly in $H_r^s(\R^3)$. The proof is completed.
\ep

\noindent
By Proposition \ref{bo5n}, for small $d\in(0,1)$ there exist $\omega_1>0, \la_2\in(0,\la_1)$ such that
\be\lab{ge}\hbox{$\|\G_\la'(u)\|_s\ge\omega_1$ for
$u\in\G_\la^{D_\la^1}\bigcap(S_1^d\setminus S_1^{\frac{d}{2}})$
and $\la\in(0,\la_2)$.}
\ee
Similar to Section 3, we have
\bo\lab{bo6n} There exists $\al_1>0$ such that for small $\la>0$,
$$
\G_\la(\g(\tau))\ge C_\la^1-\al_1\ \ \mbox{implies that}\ \ \g(\tau)\in
S_1^{\frac{d}{2}},
$$
where $\g(\tau)=U(\frac{\cdot}{\tau}), \tau\in (0,\tau_1]$. \eo

\vskip2pt
\noindent

\noindent{\bf Proof of Theorem \ref{Theorem 2} concluded.} With the help of \re{ge} and Proposition \ref{bo6n}, similarly as that in \cite{Zhang}, for $\la>0$ small enough, there exists
$\{u_n\}_n\subset \G_\la^{D_\la^1}\cap S_1^d$ such that
$\G_\la'(u_n)\rg 0$ as $n\rg\iy$. Similar as above, there exists $u_\la\in S_1^d$ with $u_\la\not\equiv0$ for small $d>0$. Moreover, up to a subsequence, $u_n\rg u_\la$ weakly in $H_r^s(\R^3)$ and a.e. in $\R^3$, and $\|u_n-u_\la\|_s\le 3d$ for $n$ large. Furthermore, $\G_\la'(u_\la)=0$.
By Lemma \ref{l3},
$$
\G_\la(u_n)=\G_\la(u_\la)+\G_\la(u_n-u_\la)+o(1).
$$
By $(H2)'$-$(H3)'$, for some $C>0$,
\begin{align*}
\G_\la(u_n-u_\la)\ge&\frac{1}{2}\int_{\R^3}(|(-\DD)^{s/2}(u_n-u_\la)|^2+\frac{1}{2}|u_n-u_\la|^2)\, \ud x\\
&-C\int_{\R^3}|u_n-u_\la|^{2^\ast_s}\, \ud x.
\end{align*}
Then by Lemma \ref{l2}, $\liminf_{n\rg\iy}\G_\la(u_n-u_\la)\ge0$ for small $d>0$. So $u_\la\in
\G_\la^{D_\la^1}\cap S_1^d$ with $\G_\la'(u_\la)=0$. Thus $u_\la$ is a nontrivial solution of \re{q2}. The asymptotic behavior of $u_\la$ follows from Proposition \ref{bo5n}. The proof is completed. \qed

\bigskip
\bigskip
\end{document}